\documentclass[a4paper,12pt]{article}
\usepackage{amsfonts}
\usepackage{amsmath}
\usepackage{latexsym}
\usepackage{graphicx}
\usepackage{amsfonts}
\usepackage{amssymb}
\usepackage{fancyhdr}
\usepackage{color}
\usepackage[backref]{hyperref}
\usepackage{pdfsync}

\providecommand{\U}[1]{\protect\rule{.1in}{.1in}}
\providecommand{\U}[1]{\protect\rule{.1in}{.1in}}
\newtheorem{theorem}{Theorem}[section]

\newtheorem{corollary}[theorem]{Corollary}

\newtheorem{definition}[theorem]{Definition}

\newtheorem{lemma}[theorem]{Lemma}

\newtheorem{remark}[theorem]{Remark}

\newcommand{\beq}{\begin{equation}}
\newcommand{\eeq}{\end{equation}}
\newcommand{\ben}{\begin{eqnarray}}
\newcommand{\een}{\end{eqnarray}}
\newcommand{\beno}{\begin{eqnarray*}}
\newcommand{\eeno}{\end{eqnarray*}}

\begin{document}
%On the uniqueness and non-uniqueness of the steady planar Navier-Stokes equations in
%an exterior domain
%Constant vorticity flows of the steady planar Navier-Stokes equations in
%an exterior domain
\title{On the uniqueness and non-uniqueness of the steady planar Navier-Stokes equations in
an exterior domain}
\author{Zhengguang Guo \thanks{%
		School of Mathematics and Statistics, Huaiyin Normal University, Huaian 223300, Jiangsu, P.R. China; E-mail: gzgmath@hytc.edu.cn}, Wendong Wang\thanks{School of  Mathematical Sciences, Dalian University of Technology, Dalian 116024, P.R. China; E-mail: wendong@dlut.edu.cn}}
\maketitle
\noindent\textbf{Abstract:} In this paper we investigate the uniqueness of solutions of the steady planar Navier-Stokes equations with different boundary conditions in the
exterior domain. For a class of incompressible flow with constant vorticity, we prove the uniqueness of the solution under the enhanced Navier boundary conditions. At the same time, some counterexamples are given to show that the uniqueness of the solution fails under the Navier boundary conditions. For the general incompressible flow with Dirichlet boundary condition, we prove various sufficient conditions for the uniqueness of the solution.

%we show some uniqueness or non-uniqueness results are obtained by given mixed slip boundary conditions. Moreover, we present some sufficient
%conditions to guarantee the triviality of steady solutions to the planar
%Navier-Stokes equations with no-slip boundary conditions in an
%exterior domain of $\mathbb{R}^{2}$.


\section{Introduction and main results}

The motion of viscous incompressible fluid past an
obstacle can be described by the Navier-Stokes equations.  More precisely, the velocity
$u$ of fluid and the pressure $\pi$ satisfy the following stationary incompressible Navier-Stokes problem:
\begin{equation}\label{EQ:NS}
\left\{\begin{array}{l}
u\cdot \nabla u+\nabla \pi=\Delta u, \quad \text{in} \,\, \Omega, \\
{\rm div }~ u=0, \hspace{2.11cm} \text{in} \,\, \Omega,
\end{array}\right.
\end{equation}
%\begin{eqnarray}
%\left( u\cdot \nabla \right) u+\nabla p&=&\Delta u\mathbf{,}  \label{EQ:NS} \\
%\text{div}u&=&0,  \label{NSE2}
%\end{eqnarray}%
where the flow domain $\Omega$ is a two-dimensional exterior domain, i.e., the complement of a bounded
domain which represents the obstacle. Without loss of generality, one assumes that $\Omega =\mathbb{R}^{2}\backslash \overline{B_1}$ and $B_{1}$ is the
disk of radius $1$ centered at the origin. There are several possibilities of
boundary conditions. According to the idea that the fluid
cannot slip on the boundary due to its viscosity, the widely used are the following Dirichlet or no-slip boundary conditions:
\begin{equation}
	u=0 \quad \text{on} \quad \partial B_1,  \label{boundary}
\end{equation}
However, in the case where the obstacles
have an approximate limit, the Dirichlet boundary conditions are no longer valid (see for example \cite{serrin}). Due to the roughness of the boundary and
the viscosity of the fluid, it is usually assumed that there is a stagnant fluid layer near the boundary, which allows the fluid to slip. This situation seems to match the reality. Then, it is really important to introduce another boundary
conditions to describe the behavior of fluid on the boundary.  In 1827, C. Navier \cite{Navier}
was the first mathematician who considered the slip phenomena and proposed the Navier-slip boundary conditions:
\begin{equation}\label{nav boundary}
	\left\{\begin{array}{l}
		u\cdot n=0,  \\
		2[D(u)\cdot n]_{\tau}+\alpha(x)u_{\tau}=0,
	\end{array}\right.
\end{equation}
where $D(u)$ is the stress tensor of fluid,  $n$ and $\tau$ are the unit outer normal vector and tangential vector of the boundary, $\alpha(x)$ is a physical parameter, which
can be a positive constant or a $L^\infty $ function on the boundary. For the far field of the fluid, one usually assumes that
\begin{equation}
u(x_1,x_2)\to \tilde{u}_\infty \quad \text{as}\quad |(x_1,x_2)|\to \infty.
\label{inf}
\end{equation}
So, it is interesting to find solutions to (\ref{EQ:NS}) with different boundary conditions. An arbitrary solution $u$ to the Navier-Stokes equations (\ref{EQ:NS}) having the finite Dirichlet integral
\begin{equation}
\int_{\Omega }|\nabla u|^{2}dx<\infty,  \label{Dirichlet}
\end{equation}%
is usually called $D$-solution \cite{galdibook}, and as is well known (see \cite{ladyzhen}), such solutions are real analytic in $\Omega$.

The study of (\ref{EQ:NS}) with conditions (\ref{boundary})(\ref{inf}) and (\ref{Dirichlet}) began with Leray \cite{Leray}
who sought the solutions as the limit of certain approximate solutions, but
the behavior of Leray solution at infinity was not found. Indeed, it was
not even apparent that Leray solution was non-trivial. The Navier-Stokes equations have been shown to have a solution by Finn and
Smith \cite{Fin} with some smallness assumptions on $|\tilde{u}_\infty|$ by
using the contraction mapping principle. Then, Amick \cite{Amick1984} proved the existence of solutions for given external forces when the exterior domain is invariant under the transformation $(x_1,x_2)\mapsto (-x_1,x_2)$, this work was generalized by Pileckas-Russo \cite{RussoMA}. Hillairet-Wittwer \cite{peterJDE} proved the existence of vanishing at infinity solutions to (\ref{EQ:NS}) with non-zero Dirichlet boundary conditions by perturbation around the radial and rotational flow $\mu x^{\perp}/|x|^2$, while the flow $\mu x^{\perp}/|x|^2$ is the exact solution decaying in the scale-critical order $O(|x|^{-1})$ under the zero flux condition. The problem of the asymptotic
behavior at infinity of an arbitrary $D$-solution $(u,\pi)$ to (\ref{EQ:NS})
 was initiated by Gilbarg-Weinberger \cite{GW} and Amick
\cite{amick}. In \cite{GW}, the authors have shown that the pressure $\pi$ has
a finite limit at infinity, and
\begin{equation*}
u(x)=o(\text{log}^{1/2}r), \quad \nabla u(x)=o(r^{-3/4}\text{log}^{9/8}r),
\end{equation*}
and
\begin{equation*}
w(x)=o(r^{-3/4}\text{log}^{1/8}r),
\end{equation*}
where $r=|x|=\sqrt{x_1^2+x_2^2}$ and $w(x)$ is the vorticity $w=\partial_2 u_1-\partial_1 u_2$.
In the elegant paper \cite{amick}, it was shown that if $(u,\pi)$ is a
solution of (\ref{EQ:NS})-(\ref{boundary}) and (\ref{Dirichlet}), then $u\in
L^{\infty }(\Omega )$ (see Theorem 12 of \cite{amick}). The assumption (\ref{boundary}) was recently removed by Korobkov-Pileckas-Russo in \cite{koro1}. Moreover, there
exists a constant vector $u_{\infty }$ such that%
\begin{equation*}
\lim_{r\rightarrow \infty }\int_{0}^{2\pi }|u\mathbf{(}r,\theta \mathbf{)}
-u_{\infty }|^{2}d\theta =0,
\end{equation*}%
and if $u_{\infty }=0,$ then $u\mathbf{(}x\mathbf{)}\rightarrow 0,$
uniformly as $|x|\rightarrow \infty ,$ where $x=(x_1,x_2)$. Some decay
properties on vorticity $w (x)$ and $\nabla u(x)$ were also obtained.
Particularly, it was proved in \cite{amick} that the following uniform limit at
infinity holds
\begin{equation*}
|u(x)|\rightarrow \left\vert u_{\infty }\right\vert \ \ \ \ \text{as}\ \ \ \
|x|\rightarrow \infty .
\end{equation*}%
Furthermore, for symmetric flow, there holds the following uniform convergence of the velocity
\begin{equation}
|u(x)\mathbf{-}u\mathbf{_{\infty }|}\rightarrow 0\ \ \ \ \text{as}\ \ \ \
|x|\rightarrow \infty .  \label{convergence}
\end{equation}%
However, the uniform convergence (\ref{convergence}) was also proved very
recently without Amick's symmetric condition or zero boundary condition on $%
\partial \Omega $ by Korobkov-Pileckas-Russo in \cite{koro}, this ensures that
the solution behaves at infinity as that of the linear Oseen equations (see,
for example, \cite{galdibook}). For asymptotic behaviour of steady  solutions to the Navier-Stokes equations at infinity, one can refer to \cite{babenko,ChoeJMFM,FarwigNNAA,SverakCPDE}. Note that the problem of the coincidence
with $u_\infty$ and the prescribed data $\tilde{u}_\infty$ is still open.
When $\tilde{u}_\infty=0$, one always has at least the trivial solution $%
u\equiv 0$, but it is not sure that:\\
{\bf Whether $u=0$ is the unique solution of (\ref%
{EQ:NS}) with the conditions of (\ref{boundary}), (\ref{inf}) and (\ref{Dirichlet})?}\\
which is exactly the conjecture raised by Amick in
\cite{amick}, and usually called the Liouville problem. There are few studies on the Liouville problem in two
dimensional exterior domains for the Navier-Stokes equations and we refer to the recent result by Korobkov-Ren in \cite{KR} for $\tilde{u}_\infty\neq 0$.

While in three
dimensional case, Galdi \cite{galdibook}
proved the Liouville type theorems by assuming that $u\in L^{9/2}(%
\mathbb{R}^{3})\cap L^{\infty }(\mathbb{R}^{3})$. Chae in \cite{ChaeCMP} showed that
the condition $\Delta u\in L^{6/5}(\mathbb{R}^{3})$ is enough to guarantee
the triviality of $u$. For related discussion, we refer to \cite{ZhangJFA,ZhangARMA,ChaeDCDS,ChaeJMAA,KorobkovJMFM,Seregin,WangJDE} and
references therein. To our knowledge, the following
Liouville problem is still open: Is a $D$-solution to (\ref{EQ:NS}) in $\mathbb{R}^{3}$, vanishing at infinity, identically zero? As a
matter of fact, this problem is also related to the very hard problem of
uniqueness of solutions to nonhomogeneous problem for the Navier-Stokes
equations. However,
the case in $\mathbb{R}^{2}$ is different, it was proved by Gilbarg-Weinberger in \cite{GW} by using the maximum
principle for the vorticity equation
\begin{equation}
\Delta w -u\cdot \nabla w =0.  \label{vort}
\end{equation}

In this paper, we focus on the Amick's conjecture for $\tilde{u}_\infty= 0$ in
\cite{amick}, and consider the uniqueness problem of steady solutions to (%
\ref{EQ:NS}) with some prescribed boundary conditions. %One one
%hand, one noted
%that the vorticity equation ({\ref{vort}}) plays an important role in
%studying the Liouville problem in $\mathbb{R}^{2}$. However, for exterior
%domains, the main difficulty is due to lack of asymptotical behavior of
%vorticity on the boundary. So we study the Liouville problem in Theorem \ref%
%{th2} for the Navier-Stokes equations (\ref{EQ:NS}) equipped
%with the following boundary conditions:
The first result is concerned with the uniqueness of the constant vorticity flow of $u_{0}=\frac{a}{2}%
(x_{2},-x_{1})$ under the following boundary condition:
\begin{equation}
u\mathbf{|}_{\partial B_{1}}=\frac{a}{2}(x_2,-x_1)\mathbf{|}_{\partial
B_{1}}, \label{BC-noslip}
\end{equation}%
where $a$ is a constant, which is stated as following for a perturbation of $L^q$ energy norm. %In this case, we perturb the Navier-Stokes
%equations around $(u_0,\pi_0)=\left(\frac{a}{2}(x_2,-x_1), \frac{a^2}{8}%
%(x_1^2+x_2^2)\right)$, then try to show the triviality of perturbed system
%using cut-off method and energy estimates.

\begin{theorem}%[Uniquenss under the no-slip boundary]%
	\label{th1} Let $(u,\pi)$ be a smooth solution of the 2D Navier-Stokes
	equations (\ref{EQ:NS}) defined on $\Omega $ and $u\in C^1(\bar{\Omega})$ satisfies the boundary
	conditions (\ref{BC-noslip}). Moreover, let $v=u-u_{0}$ with $u_{0}=\frac{a}{2}%
	(x_{2},-x_{1})$ and
	\ben\label{eq:nabla v}
	 v\in L^{q}(\Omega ),\quad 1<q\leq 2.
	\een
	Then $u\equiv u_{0}$.
\end{theorem}

The second result is concerned with the uniqueness (up to some constant) of the constant vorticity flow of $u_{0}=\frac{a}{2}%
(x_{2},-x_{1})$ under the enhanced Navier slip boundary conditions, which
is stated as following for a perturbation of $L^q$ energy norm.
Recall the Navier slip boundary conditions is as follows:
\begin{equation}
	u\cdot \vec{n}\mathbf{|}_{\partial B_{1}}=0,\text{ \ }w \mathbf{|}_{\partial B_{1}}=a,  \label{BC-slip}
\end{equation}%
and here we added an additional condition:
\ben\label{eq:enhanced slip}
\int_{\partial {B_{1}}}\frac{\partial
	v}{\partial n}\cdot vdS=0, ~{\rm or}~ \frac{\partial
	v}{\partial n} \Big| _{\partial B_1}=0, ~{\rm or} ~|v|_{\partial B_1}\equiv C.
\een
%\begin{equation}
%u\mathbf{|}_{\partial B_{1}}=\frac{a}{2}(x_2,-x_1)\mathbf{|}_{\partial
%B_{1}},\text{ \ }w \mathbf{|}_{\partial B_{1}}=a,  \label{BC-NO SLIP}
%\end{equation}%
\begin{theorem}%[Uniqueness under the enhanced Navier slip boundary]%
	\label{th2} Let $(u,\pi)$ be a smooth solution of the 2D Navier-Stokes
	equations (\ref{EQ:NS}) defined in $\Omega $ and $u\in C^2(\bar{\Omega})$ satisfies the boundary
	conditions (\ref{BC-slip}) and (\ref{eq:enhanced slip}). Moreover, let $v=u-u_{0}$ with $u_{0}=\frac{a}{2}%
	(x_{2},-x_{1})$ and
	\ben\label{eq:nabla v}
	\nabla v\in L^{q}(\Omega ),\quad 1<q<\infty.
	\een
	Then $u-u_{0}\equiv C$.
\end{theorem}

Moreover, by assuming a perturbation of $L^p$ norm under the enhanced Navier slip boundary conditions, we have the following conclusions:
\begin{theorem}%[Uniqueness under the enhanced Navier slip boundary]%
	\label{th3} Let $(u,\pi)$ be a smooth solution of the 2D Navier-Stokes
	equations (\ref{EQ:NS}) defined in $\Omega $ and $u\in C^2(\bar{\Omega})$ satisfies  the boundary
	conditions (\ref{BC-slip}) and (\ref{eq:enhanced slip}). Moreover, let $v=u-u_{0}$ with $u_{0}=\frac{a}{2}%
	(x_{2},-x_{1})$ and
	\ben\label{eq:v p}
	v\in L^{p}(\Omega ),\quad 1<p<\infty.
	\een
	Then $u-u_{0}\equiv C$.
\end{theorem}
\begin{remark}[Non-uniqueness for the Navier slip boundary condition]\label{rmk:nonunique}
	The boundary conditions in Theorem \ref{th2} and \ref{th3} seem to be sharp in a sense. Consider the usual slip condition as follows
	\beno
	u\cdot \vec{n}\mathbf{|}_{\partial B_{1}}=0,\text{ \ }w \mathbf{|}_{\partial B_{1}}=a,  \label{BC-SLIP}
	\eeno
	which is just (\ref{BC-slip}).
	At this time, one can find another non-trivial solution of (\ref{EQ:NS}), which is different from the known solution $u_{0}=\frac{a}{2}(x_{2},-x_{1})$ and $\pi_{0}=\frac{a^{2}}{8}%
	(x_{1}^{2}+x_{2}^{2})$. For example,
	\ben\label{eq:many solutions}
	u=u_0+C\frac1{x_1^2+x_2^2}(x_2,-x_1)
	\een
	solves (\ref{EQ:NS}) with the boundary conditions (\ref{BC-slip}), and $\nabla v\in L^p$ with $1<p<\infty$. However, the condition (\ref{eq:enhanced slip}) does not hold.
\end{remark}

As an immediate  corollary of Theorem \ref{th2} and Remark \ref{rmk:nonunique}, one has the following result.
\begin{corollary}%[The existence for many solutions]
	Let $(u,\pi)$ be a smooth solution of the 2D Navier-Stokes equations (\ref%
	{EQ:NS})  defined in $\Omega $ and $u_{0}=\frac{a}{2}%
	(x_{2},-x_{1})$. Moreover, $u\in C^2(\bar{\Omega})$ satisfies  the boundary
	conditions (\ref{BC-slip}). Then there exist many nontrivial solutions as in (\ref{eq:many solutions}) such that $\nabla (u-u_{0})\in
	L^{p}(\Omega )$ with $1<p<\infty $ with boundary conditions (\ref%
	{BC-slip}) and $u-u_{0}$ is vanishing at infinity.
\end{corollary}

On the other hand, we consider the special constant vorticity flow of $u_{0}=\frac{a}{2}%
(x_{2},-x_{1})$ with $a=0$ when the solution $u$ satisfies the no-slip boundary condition (\ref{boundary})
and is vanishing at infinity, i.e.,
\begin{equation}
u(x_1,x_2)\to 0 \quad \text{as}\quad |(x_1,x_2)|\to \infty.  \label{infty1}
\end{equation}
Then, some sufficient conditions which guarantee the triviality of $q$%
-generalized solutions to the Navier-Stokes equations (\ref{EQ:NS}) are established in the following theorem.
First let us recall the definition of $q$-generalized solutions.
\begin{definition}
A vector field $u:\Omega\rightarrow\mathbb{R}^{2}$ is called a $q$%
-generalized solution to (\ref{EQ:NS}), (\ref{boundary}) and (\ref{infty1}) if for some $q\in(1,\infty)$ the following properties are met:

\begin{enumerate}
\item[(i)] $u\in D_{0}^{1,q}(\Omega);$

\item[(ii)] $u$ is (weakly) divergence-free in $\Omega;$

\item[(iii)] $u$ verifies the identity%
\begin{equation*}
(\nabla u,\nabla \psi )=-(u\cdot \nabla u,\psi ),\text{ for all }\psi \in
\mathcal{D}(\Omega ).
\end{equation*}
\end{enumerate}
If $q=2,$ $u$ is usually called a generalized solution (or $D$-solution).
\end{definition}
Moreover, for the 2D Navier-Stokes equations, all $q$-generalized solutions with $q>1$ are smooth (see, for example, Chapter IX in \cite{galdibook}).

Our result is stated as follows:
\begin{theorem}%[Uniqueness under the no-slip boundary]
\label{th4} Let $(u,\pi)$ be a $q$-generalized solution to
the Navier-Stokes equations (\ref{EQ:NS}) with boundary
conditions (\ref{boundary}), (\ref{infty1}) in the exterior domain $\Omega $%
. Then, $u(x)\equiv 0$ in $\Omega $ under one of the
following conditions
\begin{enumerate}

\item $u(x)$ is a $q$-generalized solution for $1<q\leq 3/2.$

\item $u\mathbf{(}x\mathbf{)}\in BMO^{-1}(\Omega )$ for $3/2<q<2.$

\item $u\mathbf{(}x\mathbf{)}\in L^{4}(\Omega)$ for a $D$-solution.
\end{enumerate}
\end{theorem}

The rest of this paper is organized as follows. Some elementary results on
functions with finite Dirichlet energy or $q$-generalized integrals, and
the Giaquinta's iteration lemma are collected in Section \ref{sec pre},
which are important for the analysis in the rest of this paper. The proofs
of Theorems \ref{th1}, \ref{th2}, \ref{th3} and \ref{th4} are presented in Section \ref{sec-3}
-Section \ref{proof1}, respectively.

\section{Preliminaries}

\label{sec pre}

Before going to the detailed proofs of the theorems, for reader's convenience, we
would like to collect some basic lemmas, which will be used in the proof.

%Before stating the main results of this paper, the following definition and
%notations are introduced.

First, the space $BMO$ in $\Omega $ is defined as follows,
which is similar to the case in $\mathbb{R}^{2}$ as defined in \cite{H.B1}.
\begin{definition}
The space $BMO(\Omega )$ of bounded mean oscillations is the set of locally
integrable functions $f$ such that
\begin{equation}
\Vert f\Vert _{BMO}\overset{def}{=}\sup_{B}\frac{1}{|B|}\int_{B}|f-f_{B}|dx<%
\infty \quad \text{with}\quad f_{B}\overset{def}{=}\frac{1}{|B|}\int_{B}fdx.
\label{BMO defn}
\end{equation}%
The above supremum is taken over the set of Euclidean balls.
\end{definition}

It is clear that this quantity $\Vert \cdot \Vert _{BMO}$ is in general a
seminorm, unless one argues modulo constant functions, and for $f\in BMO,$
the following inequality holds true for all balls $B$%
\begin{equation}
\frac{1}{|B|}\int_{B}|f-f_{B}|^{p}dx\leq C_{p}\Vert f\Vert _{BMO}^{p}\text{ }%
,  \label{BMO pro}
\end{equation}%
where $1\leq p<\infty .$ In the following, a space that will be used is
provided by the set of functions which are derivatives of functions in $BMO$%
. More precisely, we are talking about the space introduced by Koch and
Tataru in \cite{tataru}, which is denoted by $BMO^{-1}$ (or by $\nabla BMO$)
and is defined as the space of tempered distributions $f$ such that there
exists a vector function $g=(g_{1},g_{2},g_{3})$ belonging to $BMO$ such
that $f=\nabla \cdot g$. The norm in $BMO^{-1}$ is defined by
\begin{equation*}
\Vert f\Vert _{BMO^{-1}}=\inf_{g\in BMO}\sum_{j=1}^{3}\Vert g_{j}\Vert
_{BMO}.
\end{equation*}

Second,
for $1\leq p\leq \infty ,$ let $L^{p}$ denote the usual scalar-valued and
vector-valued $L^{p}$ space over $\Omega .$ Let%
\begin{equation*}
W^{m,p}(\Omega )=\{u\in L^{p}:D^{\alpha }u\in L^{q}(\Omega ),\text{ }|\alpha
|\leq m,\text{ }m\in \mathbb{N}\}.
\end{equation*}%
When $p=2,$ one abbreviates $H^{m}(\Omega )=W^{m,2}(\Omega ).$ If $q\in
\lbrack 1,n),$ the space $D_{0}^{1,q}(\Omega )$ is the following:%
\begin{equation*}
D_{0}^{1,q}(\Omega )=\{u\in D^{1,q}(\Omega ):\left\Vert u\right\Vert
_{nq/(n-q)}<\infty ,\text{ }\varphi u\in W_{0}^{1,q}(\Omega ),\text{ for all
}\varphi \in C_{0}^{\infty }(\mathbb{R}^{2})\},
\end{equation*}%
if $q\geq n,$ and complementary set $\Omega ^{c}\supset B_{a},$ for some $%
a>0:$%
\begin{equation*}
D_{0}^{1,q}(\Omega )=\{u\in D^{1,q}(\Omega ):\varphi u\in W_{0}^{1,q}(\Omega
),\text{ for all }\varphi \in C_{0}^{\infty }(\mathbb{R}^{2})\},
\end{equation*}%
where $W_{0}^{1,q}(\Omega )$ is the closure of $C_{0}^{\infty }(\Omega )$ in
the classical Sobolev space $W^{1,q}(\Omega )$, and%
\begin{align*}
D^{1,q}(\Omega )& =\{u\in L_{loc}^{1}(\Omega ):\nabla u\in L^{q}(\Omega )\},
\\
\mathcal{D}(\Omega )& =\{\psi \in C_{0}^{\infty }(\Omega ):\text{div}u=0%
\text{ in }\Omega \mathcal{\}}.
\end{align*}

Finally, let us  recall some necessary lemmas. The first one is a lemma of Gilbarg-Weinberger in \cite{GW} about the
decay of functions with finite Dirichlet integrals.
\begin{lemma}[Lemma 2.1, \protect\cite{GW}]
\label{GW} Let a $C^{1}$ vector-valued function $f(x)=(f_{1},f_{2})(x)=f(r,%
\theta )$ with $r=|x|$ and $x_{1}=r\cos \theta $. There holds finite
Dirichlet integral in the range $r>r_{0}$, that is
\begin{equation*}
\int_{r>r_{0}}|\nabla f|^{2}\,dx<\infty .
\end{equation*}%
Then, we have
\begin{equation*}
\lim_{r\rightarrow \infty }\frac{1}{\ln r}\int_{0}^{2\pi }|f(r,\theta
)|^{2}d\theta =0.
\end{equation*}
\end{lemma}

For general energy integrals, we have the following:
\begin{lemma}[Theorem II.9.1, \protect\cite{galdibook}]
\label{Galdi lemma} Let $\Omega \subset \mathbb{R}^{2}$ be an exterior
domain.

(i) Let
\begin{equation*}
\nabla f\in L^{r}\cap L^{p}(\Omega ),
\end{equation*}%
for some $1\leq r<2<p<\infty $. Then there exists $f_{0}\in \mathbb{R}$ such
that
\begin{equation*}
\lim_{|x|\rightarrow \infty }|f(x)-f_{0}|=0,
\end{equation*}%
uniformly.\newline

(ii) Let
\begin{equation*}
\nabla f\in L^{2}\cap L^{p}(\Omega ),
\end{equation*}%
for some $2<p<\infty $. Then
\begin{equation*}
\lim_{|x|\rightarrow \infty }\frac{|f(x)|}{\sqrt{\ln (|x|)}}=0,
\end{equation*}%
uniformly.\newline

(iii) Let
\begin{equation*}
\nabla f\in L^{p}(\Omega ),
\end{equation*}%
for some $2<p<\infty $. Then
\begin{equation*}
\lim_{|x|\rightarrow \infty }\frac{|f(x)|}{|x|^{\frac{p-2}{p}}}=0,
\end{equation*}%
uniformly.
\end{lemma}

Moreover, we recall a lemma from \cite{Wang2}.
\begin{lemma}
\label{GW2} Let a $C^{1}$ vector-valued function $f(x)=(f_{1},f_{2})(x)=f(r,%
\theta )$ with $r=|x|$ and $x_{1}=r\cos \theta $. There holds
\begin{equation*}
\int_{r>r_{0}}|\nabla f|^{q}dx<\infty ,\quad 1<q<2.
\end{equation*}%
Then, we have
\begin{equation*}
\limsup_{r\rightarrow \infty }\int_{0}^{2\pi }|f(r,\theta )|^{q}d\theta
<\infty .
\end{equation*}
\end{lemma}

The following one is the Giaquinta's iteration lemma, which gives the
estimates of the $L^{2\text{ }}$norm of $\nabla u$ in our proof.
\begin{lemma}[Lemma 3.1, Page 161, \protect\cite{G83}]
\label{Gia iteration} Let $f(t)$ be a nonnegative bounded function defined
in $[r_{0},r_{1}],$ $r_{0}\geq 0.$ Suppose that for $r_{0}\leq t<s\leq r_{1}$
we have%
\begin{equation*}
f(t)\leq \lbrack A(s-t)^{-\alpha }+B]+\theta f(s),
\end{equation*}%
where $A,$ $B,$ $\alpha ,$ $\theta $ are nonnegative constants with $0\leq
\theta <1.$ Then for all $r_{0}\leq \rho <R\leq r_{1}$ we have%
\begin{equation*}
f(\rho )\leq c[A(R-\rho )^{-\alpha }+B],
\end{equation*}%
where $c$ is a constant depending on $\alpha $ and $\theta .$
\end{lemma}

Next, the following Gagliardo-Nirenberg inequality (see \cite[Theorem 10.1,
Page 27]{Friedman}) will be frequently used.
\begin{lemma}
\label{GN1} Let $\Omega _{0}\subset \mathbb{R}^{2}$ be a bounded smooth
domain. Assume that $1\leq q,r\leq \infty ,$ and $j,m$ are arbitrary
integers satisfying $0\leq j<m.$ If $v\in W^{m,r}(\Omega _{0})\cap
L^{q}(\Omega _{0}),$ then we have%
\begin{equation*}
\left\Vert D^{j}v\right\Vert _{L^{p}}\leq C\left\Vert v\right\Vert
_{L^{q}}^{1-a}\left\Vert v\right\Vert _{W^{m,r}}^{a},
\end{equation*}%
where%
\begin{equation*}
-j+\frac{2}{p}=(1-a)\frac{2}{q}+a\left( -m+\frac{2}{r}\right) ,
\end{equation*}%
and%
\begin{equation*}
a\in \left\{
\begin{array}{l}
\left[ \frac{j}{m},1\right) ,\text{ if }m-j-\frac{2}{r}\text{ is a
nonnegative integer,} \\
\left[ \frac{j}{m},1\right] ,\text{ otherwise,}%
\end{array}%
\right.
\end{equation*}%
the constant $C$ depends only on $m,j,q,r,a,$ and $\Omega _{0}.$
\end{lemma}

\section{Proof of Theorem \ref{th1}} \label{sec-3}
\smallskip

%\textbf{Proof of Case IV}

%In this case, we want to prove that $\nabla v\in L^{2}(\Omega ),$ then the
%proof is complete by \textbf{Case II}.

%In this case, we want to prove that $\nabla v\in L^2(\Omega),$ then the proof is complete by Theorem \ref{th2}.

%{\bf Proof of Theorem \ref{th3}.}

We perturb the Navier-Stokes equations around $(u_0,\pi_0)=\left(\frac{a}{2}(x_2,-x_1), \frac{a^2}{8}%
(x_1^2+x_2^2)\right)$, then try to show the triviality of the perturbed system. Since the constant vorticity flow of $(u_{0},\pi_{0})$ solves the system (\ref{EQ:NS}),  then $%
v=u-u_{0}$ and $\pi_{1}=\pi-\pi_{0}$ satisfy the following system
\begin{equation}
	\left\{
	\begin{array}{ll}
		-\Delta v+u\cdot \nabla v+v\cdot \nabla u_{0}+\nabla \pi_{1}=0, \\
		\mathrm{div}~v=0,%
	\end{array}%
	\right.  \label{NS-2D-r}
\end{equation}%
with the boundary condition
\ben\label{eq:BC-v}
v|_{\partial B_1}=0,
\een
due to (\ref{BC-noslip}). Next we show that $v\equiv0$ under the assumptions of Theorem \ref{th1}.

%with the vorticity $\tilde{w}\doteq\partial _{2}v_{1}-\partial _{1}v_{2}=w-a$.
%Then the equation of the vorticity $\tilde{w}$ is as follows:
%\begin{equation}
%	-\Delta \tilde{w}+v\cdot \nabla \tilde{w}+u_{0}\cdot \nabla \tilde{w}=0.
%	\label{NSv-2D}
%\end{equation}
%Furthermore, let $v'=\psi v$, where $\psi$ a smooth cut-off function with $0\leq \psi \leq
%1$ satisfying
%\begin{equation}
%	\psi(x)=\psi(|x|)=\left\{
%	\begin{array}{l}
%		0,\quad |x|\leq 2, \\
%		1, \quad |x|\geq 3.%
%	\end{array}%
%	\right.  \label{cutoff psi}
%\end{equation}
%Then $v'\in C^\infty(\mathbb{R}^2)$ and $v'(x)\equiv v(x)$ for $|x|\geq 3.$ Similarly, define the vorticity $w'=\partial _{2}v'_{1}-\partial _{1}v'_{2}$.

First, we introduce a cut-off function $\phi(x)\in C_0^\infty(B_R)$ with $0\leq \phi\leq 1$ satisfying the following two properties:
\begin{enumerate}
\item[i).]
$\phi$ is radially decreasing and satisfies
\ben\label{eq:phi} \phi(x)=\phi(|x|)=\left\{
\begin{aligned}
&1,\quad |x|\leq \rho,\\
&0, \quad |x|\geq\tau,
\end{aligned}
\right. \een
where $1<\frac{R}{2}\leq \frac23 \tau\leq  \rho<\tau\leq R$;
\item[ii).]
 $|\nabla\phi(x)|\leq \frac{C}{\tau-\rho}$, $|\nabla^2\phi(x)|\leq \frac{C}{(\tau-\rho)^2}$ for all $x\in \mathbb{R}^2$.
\end{enumerate}

Second, due to the choosing of $\phi$ and (\ref{eq:BC-v}), one has
\beno
\int_{B_\tau\setminus B_1}\nabla \cdot(\phi v)dx=-\int_{\partial B_1}n\cdot v\phi dS=-\phi(1)\int_{\partial B_1}n\cdot v dS=0.
\eeno

We recall now the Bogovski\u{i} problem:
\begin{equation}
	\nabla\cdot \hat{w}=\nabla\cdot[\phi {v}]. \label{bogovski}
\end{equation}%
where a vector-valued function $\hat{w}: B_\tau\setminus{B_{\frac23\tau}}\rightarrow \mathbb{R}^2$.
Due to Bogovski\u{i}'s result in \cite{bogov} (see also, Theorem III 3.1 in \cite{galdibook}), there exists a constant $C(s)$ and a vector-valued function $\hat{w}$ such that
$\hat{w}\in W^{1,s}_0(B_\tau\setminus{B_{\frac23\tau}})$ and
(\ref{bogovski}) holds. Furthermore, we obtain
\begin{align}
	\label{esti-ws'}
	\int_{B_\tau\setminus{B_{\frac23\tau}}}|\nabla \hat{w}|^s\,dx
	\leq C(s)\int_{B_\tau}|\nabla\phi\cdot { v}|^s\,dx.
\end{align}
Making the inner products $(\phi {v}-\hat{w})$ on both sides of the equation (\ref{NS-2D-r}), by $\nabla\cdot \hat{w}=\nabla\cdot[\phi{v}]$ we have
\begin{align*}
&\hspace{-2mm}\int_{B_\tau\setminus{\overline{B_1}}}\phi|\nabla v|^2\,dx
\\&= -\int_{B_\tau\setminus{\overline{B_1}}}\nabla\phi\cdot\nabla v \cdot {v} \,dx+\int_{B_\tau\setminus{\overline{B_1}}}\nabla \hat{w}:\nabla v  \,dx
-\int_{B_\tau\setminus{\overline{B_1}}}u\cdot\nabla v \cdot \phi {v} \,dx\\
& \hspace{5mm}+\int_{B_\tau\setminus{\overline{B_1}}}u\cdot\nabla v \cdot \hat{w} \,dx-\int_{B_\tau\setminus{\overline{B_1}}}v\cdot\nabla u_0 \cdot \phi {v} \,dx+\int_{B_\tau\setminus{\overline{B_1}}}v\cdot\nabla u_0 \cdot \hat{w} \,dx \\
&\doteq  I_1+\cdots+I_6,
\end{align*}
For the term $I_1$, it follows from H\"{o}lder's inequality that
\beno
|I_1|\leq \frac{C}{\tau-\rho}\left(\int_{B_\tau\setminus{\overline{B_1}}}|\nabla v|^2\,dx\right)^{\frac12}\left(\int_{B_\tau\setminus{B_{\frac23\tau}}}| {v}|^2\,dx\right)^{\frac12}.
\eeno
For the term $I_2$, H\"{o}lder's inequality and (\ref{esti-ws'}) imply that
\begin{align*}
|I_2|&\leq C\left(\int_{B_\tau\setminus{\overline{B_1}}}|\nabla v|^2\,dx\right)^{\frac12}\|\nabla \hat{w}\|_{L^{2}(B_\tau\setminus{\overline{B_1}})}\\
&\leq  \frac{C}{\tau-\rho}\|\nabla v\|_{L^2(B_\tau\setminus{\overline{B_1}})} \|{v}\|_{L^2(B_\tau\setminus{B_{\frac23\tau}})}.
\end{align*}
By integration by parts and (\ref{esti-ws'}), we find that
\beno
|I_3|&=&\left|\int_{B_\tau\setminus{\overline{B_1}}}u\cdot\nabla v \cdot \phi {v} \,dx\right|
=\left|\int_{B_\tau\setminus{\overline{B_1}}}v\cdot\nabla v \cdot \phi {v} \,dx\right|\\
&\leq& \frac{C}{\tau-\rho}\|v\|^{3}_{L^{3}(B_\tau\setminus{B_{\frac23\tau}})},
\eeno
and
\beno
|I_4|&\leq& \frac{C}{\tau-\rho}\|v\|^{3}_{L^{3}(B_\tau\setminus{B_{\frac23\tau}})}+\left|\int_{B_\tau\setminus{\overline{B_1}}}u_0\cdot\nabla \hat{w} \cdot v dx\right|\\
&\leq & \frac{C}{\tau-\rho}\|v\|^{3}_{L^{3}(B_\tau\setminus{B_{\frac23\tau}})} +C\frac{R}{\tau-\rho}\left(\int_{B_\tau\setminus{B_{\frac23\tau}}}| {v}|^2\,dx\right).
\eeno
Moreover, $I_5=0$ due to the antisymmetric matrix $\nabla u_0$, and
\beno
|I_6|&\leq& \left|\int_{B_\tau\setminus{\overline{B_1}}}v\cdot\nabla \hat{w} \cdot u_0 dx\right|
\leq C\frac{R}{\tau-\rho}\left(\int_{B_\tau\setminus{B_{\frac23\tau}}}| {v}|^2\,dx\right)
\eeno

Combining the estimates of $I_1$--$I_6$, we get
\begin{align}\label{eq:energy}
&\int_{B_\tau\setminus{\overline{B_1}}}\phi|\nabla v|^2\,dx \notag \\
&\leq\frac14\|\nabla v\|_{L^2(B_\tau\setminus{\overline{B_1}})}^2+ \frac{CR^2}{(\tau-\rho)^{2}}\|v\|^{2}_{L^{2}(B_\tau\setminus{B_{\frac23\tau}})}+\frac{C}{\tau-\rho}\|v\|^{3}_{L^{3}(B_\tau\setminus{B_{\frac23\tau}})}.
\end{align}

%\bigskip
%{\bf \underline{Case of $3\leq p\leq 6$}}.
%\bigskip
%
%Applying Lemma \ref{Gia iteration} to (\ref{eq:energy}), we have
%\ben\label{eq:L2 estimate}
%\int_{B_\rho\backslash{\overline{B_1}}}|\nabla v|^2\,dx\leq \frac{C}{(\tau-\rho)^{2}}\|v\|^{2}_{L^{2}(B_R\setminus{B_{R/2}})}+ \frac{C}{\tau-\rho}\|v\|^{3}_{L^{3}(B_R\setminus{B_{R/2}})}+C,
%\een
%which implies that $v\equiv 0$ due to Theorem \ref{th2}.

Finally, we deal with
the case of $v\in L^p$ with  $1< p\leq 2$. Recall that the following Poincar\'{e}-Sobolev inequality holds (see, for example, Lemma \ref{GN1} or Theorem 8.11 and 8.12 \cite{LL})
\beno
\label{eq:poincare-sobolev}
\|f\|_{L^2(B_\tau\setminus{B_{\frac23\tau}})}\leq C \|\nabla f\|_{L^2(B_\tau\setminus{B_{\frac23\tau}})}^{1-\frac{p}2}\|f\|_{L^p(B_\tau\setminus{B_{\frac23\tau}})}^{\frac{p}2}+C\tau^{1-\frac2{p}}\|f\|_{L^p(B_\tau\setminus{B_{\frac23\tau}})},
\eeno
and
\beno
\label{eq:poincare-sobolev}
\|f\|_{L^3(B_\tau\setminus{B_{\frac23\tau}})}\leq C \|\nabla f\|_{L^2(B_\tau\setminus{B_{\frac23\tau}})}^{1-\frac{p}3}\|f\|_{L^p(B_\tau\setminus{B_{\frac23\tau}})}^{\frac{p}3}+C\tau^{\frac23-\frac2{p}}\|f\|_{L^p(B_\tau\setminus{B_{\frac23\tau}})},
\eeno
which imply that
\ben\label{eq:energy-'}
&&\int_{B_\tau\setminus{\overline{B_1}}}\phi|\nabla v|^2\,dx\nonumber\\
&\leq&\frac12\|\nabla v\|_{L^2(B_\tau\setminus{\overline{B_1}})}^2+ \frac{CR^2}{(\tau-\rho)^{2}}\|v\|^{2}_{L^{2}(B_\tau\setminus{B_{\frac23\tau}})}+\frac{C}{\tau-\rho}\|v\|^{3}_{L^{3}(B_\tau\setminus{B_{\frac23\tau}})}\nonumber\\
&\leq &\frac12\|\nabla v\|_{L^2(B_\tau\setminus{\overline{B_1}})}^2+C \frac{R^2}{(\tau-\rho)^{2}}\left(\|\nabla v\|_{L^2(B_\tau\setminus{B_{\frac23\tau}})}^{2-p}\|v\|_{L^p(B_\tau\setminus{B_{\frac23\tau}})}^{p}+C\tau^{2-\frac{4}{p}}\|v\|_{L^p(B_\tau\setminus{B_{\frac23\tau}})}^2
 \right)\nonumber\\
 &&+\frac{C}{\tau-\rho}\left(  \|\nabla v\|_{L^2(B_\tau\setminus{B_{\frac23\tau}})}^{3-p}\|v\|_{L^p(B_\tau\setminus{B_{\frac23\tau}})}^{p}+\tau^{2-\frac{6}{p}}\|v\|_{L^p(B_\tau\setminus{B_{\frac23\tau}})}^3 \right).
\een
It follows from Young's inequality and $v\in L^p$ for $1<p\leq 2$ that
\beno\label{eq:energy-}
&&\int_{B_\tau\setminus{\overline{B_1}}}\phi|\nabla v|^2\,dx\\
&\leq &\frac34\|\nabla v\|_{L^2(B_\tau\setminus{\overline{B_1}})}^2+C \left(\frac{R^2}{(\tau-\rho)^{2}}\right)^{\frac2p}+C\frac{R^{4-\frac{4}{p}}}{(\tau-\rho)^{2}}+\frac{C}{(\tau-\rho)^{\frac{2}{p-1}}}+\frac{C\tau^{2-\frac{6}{p}}}{\tau-\rho}.
\eeno
Applying Lemma \ref{Gia iteration}, we have
\beno\label{eq:L2 estimate}
\int_{B_{R/2}\backslash{\overline{B_1}}}|\nabla v|^2\,dx\leq C,~~{\rm for~any~} R>2.
\eeno
Using this and (\ref{eq:energy-'}) again, by taking $\tau=2\rho=R\rightarrow\infty$ we have
\ben\label{eq:L2 estimate}
\int_{\mathbb{R}^2\backslash{\overline{B_1}}}|\nabla v|^2\,dx=0,
\een
which implies $v\equiv 0$ due to (\ref{BC-noslip}). Thus the proof is complete.

\section{Proof of Theorem \protect\ref{th2}}

\label{proof2}
Recall that $%
v=u-u_{0}$ and $\pi_{1}=\pi-\pi_{0}$ satisfy (\ref{NS-2D-r}).
Define the vorticity $\tilde{w}\doteq\partial _{2}v_{1}-\partial _{1}v_{2}=w-a$.
Then the equation of the vorticity $\tilde{w}$ is as follows:
\begin{equation}
	-\Delta \tilde{w}+v\cdot \nabla \tilde{w}+u_{0}\cdot \nabla \tilde{w}=0.
	\label{NSv-2D}
\end{equation}
Furthermore, let $v'=\psi v$, where $\psi$ a smooth cut-off function with $0\leq \psi \leq
1$ satisfying
\begin{equation}
	\psi(x)=\psi(|x|)=\left\{
	\begin{array}{l}
		0,\quad |x|\leq 2, \\
		1, \quad |x|\geq 3.%
	\end{array}%
	\right.  \label{cutoff psi}
\end{equation}
Then $v'\in C^\infty(\mathbb{R}^2)$ and $v'(x)\equiv v(x)$ for $|x|\geq 3.$ Similarly, define the vorticity $w'=\partial _{2}v'_{1}-\partial _{1}v'_{2}$, then $w'(x)\equiv \tilde{w}(x)$ for $|x|\geq 3.$

\bigskip
\underline{\textbf{Step 1. Case of $2<q<\infty .$}}
\bigskip

Let $\eta (x_1,x_2)\in
C_{0}^{\infty }(\mathbb{R}^{2})$ be a cut-off function with $0\leq \eta \leq
1$ satisfying $\eta (x)=\eta _{1}(\frac{|x|}{R})$, where
\begin{equation}
	\eta _{1}(t)=\left\{
	\begin{array}{l}
		1,\quad |t|\leq 1, \\
		0, \quad |t|>2.%
	\end{array}%
	\right.  \label{cutoff eta}
\end{equation}
Multiply $q\eta |w-a|^{q-2}(w-a)$ on both sides of (\ref{NSv-2D}), then we
have
\begin{align}
	& \frac{4(q-1)}{q}\int_{\Omega }|\nabla (|w-a|^{\frac{q}{2}})|^{2}\eta dx
	\notag  \label{eq:energy estimate-w} \\
	& \leq \int_{\Omega }|w-a|^{q}\Delta \eta dx+\int_{\Omega
	}|w-a|^{q}v\cdot \nabla \eta dx  \notag \\
	& +\int_{\Omega }|w-a|^{q}u_{0}\cdot \nabla \eta dx\triangleq
	K_{1}+K_{2}+K_{3}.
\end{align}%
Since $\tilde{w}=w-a\in L^{q}$ by (\ref{eq:nabla v}), obviously $K_{1}\rightarrow 0$ as $R\rightarrow \infty $. For the term $K_{2}$%
, due to (iii) in Lemma \ref{Galdi lemma} and (\ref{eq:nabla v}), for large $R>0$
we have
\begin{equation*}
	|v(x_1,x_2)|\leq |(x_1,x_2)|^{1-\frac{2}{q}}.
\end{equation*}%
Thus we have
\begin{equation*}
	K_{2}\leq CR^{(1-\frac{2}{q})-1}\rightarrow 0,
\end{equation*}%
as $R\rightarrow \infty $. It is worth noting that the third term is
vanishing, since $u_{0}=\frac{a}{2}(x_2,-x_1)$ belongs to the tangent vector
and $\nabla \eta $ is the radial vector.
Consequently, we get $\nabla (|\tilde{w}|^{\frac{q}{2}})\equiv 0,$ which
implies that $\tilde{w}\equiv 0$ by (\ref{BC-slip}). Due to $\mathrm{div}~v=0$, it follows that
\begin{equation*}
	\Delta v\equiv 0,\quad \mathrm{in}~\Omega .
\end{equation*}%
\textbf{Claim that}:
\begin{equation}
	v\equiv C, \quad \mathrm{in}~\Omega. \label{v vanishing}
\end{equation}

Firstly, due to $\Delta v'=\nabla ^{\perp }w'$, there holds
\ben\label{eq:w'}
	\Vert \nabla v'\Vert _{L^{\frac32}(\mathbb{R}^{2})}+\Vert \nabla ^{2}v'%
	\Vert _{L^{p}(\mathbb{R}^{2})}\leq C(\Vert w'\Vert _{L^{\frac32}(\mathbb{R}%
		^{2})}+\Vert \nabla w'\Vert _{L^{p}(\mathbb{R}^{2})})<\infty ,
\een
by the help of Calder\'{o}n-Zygmund estimates, since $w'=\tilde{w}\equiv 0 $ for $|x|\geq 3$. Due to Lemma \ref{Galdi lemma}%
, there exists a constant vector $v_{0}$ such that
\begin{equation*}
	\lim_{|x|\rightarrow \infty }|v'-v_{0}|=0,
\end{equation*}%
uniformly, which implies that
\begin{equation}
	\Vert v\Vert _{L^{\infty }(\Omega)}\leq C_{0}.  \label{infty bound}
\end{equation}

Secondly, for any $r>1$, by (\ref{eq:enhanced slip}) we have
\begin{align*}
	0& =\int_{B_{r}\setminus {B_{1}}}\Delta v\cdot vdx \\
	& =-\int_{B_{r}\setminus {B_{1}}}|\nabla v|^{2}dx+\int_{\partial B_{r}}%
	\frac{\partial v}{\partial n}\cdot vdS-\int_{\partial {B_{1}}}\frac{\partial
		v}{\partial n}\cdot vdS \\
	& =-\int_{B_{r}\setminus {B_{1}}}|\nabla v|^{2}dx+\frac{r}{2}%
	\int_{\partial B_{1}}\frac{\partial }{\partial r}\left[|v(rz)|^{2}\right]dS_{z},
\end{align*}%
which yields that
\begin{equation*}
	rG^{\prime }(r)=2\int_{B_{r}\setminus {B_{1}}}|\nabla v|^{2}dx,
\end{equation*}%
provided that
\begin{equation*}
	\int_{\partial B_{1}}\left[|v(rz)|^{2}\right]dS_{z}=G(r).
\end{equation*}%
Then by solving the ODE equation we have
\begin{equation*}
	G(r)\geq G(r_0)+\left( 2\int_{B_{r_{0}}\setminus {B_{1}}}|\nabla
	v|^{2}dxdy\right) \ln {\frac{r}{r_0}}
\end{equation*}%
for any $r>r_{0}>1.$ Note that (\ref{infty bound}) implies that $G(r)\leq C_0^2|\partial B_1|
$ for any all $r>1$, then
\begin{equation*}
	\int_{B_{r}\setminus {B_{1}}}|\nabla v|^{2}dxdy=0.
\end{equation*}
That is to say
\begin{equation*}
	\nabla v\equiv 0, \quad \text{in} \quad \Omega ,
\end{equation*}%
and $v\equiv v_{0}$. Thus we have $%
v\equiv C.$ The proof of (\ref{v vanishing}) is complete.
%Hence $v$ and $%p_{1}$ are constants.

\bigskip

\underline{\textbf{Step 2. Case of $1<q\leq 2.$}}

\bigskip

We take a cut-off function $\phi $ as in (\ref{eq:phi}).
%\begin{enumerate}
%	\item[i).] Let $r=\sqrt{x_1^{2}+x_2^{2}}$. $\phi $ is radially decreasing
%	and satisfies
%	\begin{equation*}
%		\phi (x_1,x_2)=\phi (r)=\left\{
%		\begin{array}{c}
%			1,\text{ \ }r\leq \rho , \\
%			0,\text{ \ }r\geq \tau ,%
%		\end{array}%
%		\right.
%	\end{equation*}%
%	where $1<\frac{R}{2}\leq \frac{2}{3}\tau <\frac{3}{4}R\leq \rho <\tau \leq R$%
%	;
%	
%	\item[ii).] $|\nabla\phi(x)|\leq\frac{C}{\tau-\rho}$ for all $(x_1,x_2)\in
%	\mathbb{R}^{2}$.
%\end{enumerate}

Multiplying both sides of (\ref{NSv-2D}) by $\phi (w-a)$ and then applying
integration by parts, we arrive at
\begin{align}
	& \hspace{-1mm}\int_{B_{1}^{c}}\phi |\nabla w|^{2}\,dx  \notag
	\label{eq:vorticity equ} \\
	& =-\int_{B_{1}^{c}}\nabla w\cdot \nabla \phi (w-a)\,dx+\frac{1}{2}%
	\int_{B_{1}^{c}}v\cdot \nabla \phi (w-a)^{2}\,dx+\frac{1}{2}%
	\int_{B_{1}^{c}}u_{0}\cdot \nabla \phi (w-a)^{2}\,dx  \notag \\
	& \doteq I_{1}^{\prime }+I_{2}^{\prime }+I_{3}^{\prime }.
\end{align}%
In what follows we shall estimate $I_{j}^{\prime }$ for $j=1,2,3$ one by
one. As in Step 1, $I_{3}^{\prime }=0.$

For the term $I_{1}^{\prime }$, by H\"{o}lder's inequality we have
\begin{equation*}
	I_{1}^{\prime }\leq \frac{C}{\tau -\rho }\Vert \nabla w\Vert _{L^{2}(B_{\tau
		}\setminus{B_1})}\Vert w-a\Vert _{L^{2}(B_{\tau }\setminus {B_{\frac{2}{3}\tau }})},
\end{equation*}%
Using the following multiplicative Gagliardo-Nirenberg inequality again
\ben\label{eq:L2 interpolation}
\Vert w-a\Vert _{L^{2}(B_{\tau }\setminus {B_{\frac{2}{3}\tau }})}\leq
C\Vert \nabla w\Vert _{L^{2}(B_{\tau }\setminus {B_{\frac{2}{3}\tau }})}^{1-%
	\frac{q}{2}}\Vert w-a\Vert _{L^{q}(B_{\tau }\setminus {B_{\frac{2}{3}\tau }}%
	)}^{\frac{q}{2}}+C\tau ^{1-\frac{2}{q}}\Vert w-a\Vert _{L^{q}(B_{\tau
	}\setminus {B_{\frac{2}{3}\tau }})},\nonumber\\
\een
which yields that
\begin{equation}
	I_{1}^{\prime }\leq \frac{1}{8}\int_{B_{\tau
		}\setminus{B_1}}|\nabla w|^{2}\,dx+\frac{C}{%
		(\tau -\rho )^{\frac{4}{q}}}+\frac{C\tau ^{2-\frac{4}{q}}}{(\tau -\rho )^{2}}%
	,  \label{estimate of I1}
\end{equation}%
by noting that $\Vert w-a\Vert _{L^{q}(B_{\tau }\setminus {B_{\frac{2}{3}%
			\tau }})}<\infty .$

For the terms $I_{2}^{\prime }$, let
\begin{equation*}
	\bar{f}(r)=\frac{1}{2\pi }\int_{0}^{2\pi }f(r,\theta )d\theta ,
\end{equation*}%
then by Writinger's inequality (for example, for $p=2$ see Chapter II.5 \cite%
{galdibook}) we have
\begin{equation}
	\int_{0}^{2\pi }|f-\bar{f}|^{p}\,d\theta \leq C(p)\int_{0}^{2\pi }|\partial
	_{\theta }f|^{p}d\theta ,  \label{Wirtinger}
\end{equation}%
for $1\leq p<\infty $.

Then by using (\ref{Wirtinger}), Lemma \ref{GW} and Lemma \ref{GW2} we have
\begin{align*}
	I_{2}^{\prime }& \leq \left\vert \int_{\mathbb{R}^{2}}(w-a)^{2}\,(v-\bar{v}%
	)\cdot \nabla \phi \,dx\right\vert +\left\vert \int_{\mathbb{R}%
		^{2}}(w-a)^{2}\,\bar{v}\cdot \nabla \phi \,dx\right\vert \\
	& \leq \frac{C}{\tau -\rho }\left( \int_{B_{\tau }\setminus {B_{\frac{2}{3}%
				\tau }}}(w-a)^{2q^{\prime }}\right) ^{\frac{1}{q^{\prime }}}\left( \int_{%
		\frac{2}{3}\tau <r<\tau }\int_{0}^{2\pi }|v(r,\theta )-\bar{v}|^{q}\,d\theta
	\,rdr\right) ^{\frac{1}{q}} \\
	& +\frac{C}{\tau -\rho }\int_{B_{\tau }\setminus {B_{\frac{2}{3}\tau }}%
	}(w-a)^{2}\left( \int_{0}^{2\pi }|v(r,\theta )|^{q}\,d\theta \right) ^{\frac{%
			1}{q}}\,dx \\
	& \leq \frac{CR}{\tau -\rho }\left( \int_{B_{\tau }\setminus {B_{\frac{2}{3}%
				\tau }}}(w-a)^{2q^{\prime }}\right) ^{\frac{1}{q^{\prime }}}\left( \int_{%
		\frac{2}{3}\tau <r<\tau }\frac{1}{r^{q}}\int_{0}^{2\pi }|\partial _{\theta
	}v|^{q}d\theta \,rdr\right) ^{\frac{1}{q}} \\
	& +C\frac{(\ln R)^{\frac{1}{2}}}{\tau -\rho }\int_{B_{\tau }\setminus {B_{%
				\frac{2}{3}\tau }}}(w-a)^{2}\,dx.
\end{align*}%
Using Gagliardo-Nirenberg inequality again, one has%
\begin{equation}\label{Lqinterpolation}
\Vert w-a\Vert _{L^{2q^{\prime }}(B_{\tau }\setminus {B_{\frac{2}{3}\tau }}%
	)}\leq C\Vert \nabla w\Vert _{L^{2}(B_{\tau }\setminus {B_{\frac{2}{3}\tau }}%
	)}^{1-\frac{q}{2q^{\prime }}}\Vert w-a\Vert _{L^{q}(B_{\tau }\setminus {B_{%
			\frac{2}{3}\tau }})}^{\frac{q}{2q^{\prime }}}+C\tau ^{1-\frac{3}{q}}\Vert
w-a\Vert _{L^{q}(B_{\tau }\setminus {B_{\frac{2}{3}\tau }})}.
\end{equation}
It follows from (\ref{eq:L2 interpolation}) and (\ref{Lqinterpolation}) that
\begin{align}
	I_{2}^{\prime }& \leq \frac{1}{8}\left( \int_{B_{\tau
		}\setminus{B_1}}|\nabla
	w|^{2}\right) +C\left( \frac{R}{\tau -\rho }\right) ^{\frac{%
			2q^{\prime }}{q}}\left(\Vert \nabla v\Vert
	_{L^{q}(B_{\tau }\setminus {B_{\frac{2}{3}\tau }})}\right) ^{\frac{%
			2q^{\prime }}{q}+2}+CR^{3-\frac{6}{q}} ({\tau -\rho })^{-1} \notag \\
	& +C\left( \frac{\sqrt{\ln R}}{\tau -\rho }\right) ^{\frac{2}{q}}+C\left(
	\frac{\sqrt{\ln R}}{\tau -\rho }\right) \tau ^{2-\frac{4}{q}},
	\label{estimate of I2}
\end{align}%
where we used the boundedness of $\Vert \nabla v\Vert _{L^{q}(B_{\tau }\setminus {%
		B_{\frac{2}{3}\tau }})}$.

Collecting the estimates of $I_{1}^{\prime },I_{2}^{\prime }$, by (\ref%
{estimate of I1}) and (\ref{estimate of I2}) we have
\begin{align*}
	& \int_{B_{\rho }\setminus{B_1}}|\nabla w|^{2}dx \\
	& \leq \frac{1}{2}\int_{B_{\tau
		}\setminus{B_1}}|\nabla w|^{2}+\frac{C}{(\tau -\rho )^{%
			\frac{4}{q}}}+\frac{C\tau ^{2-\frac{4}{q}}}{(\tau -\rho )^{2}}+CR^{3-\frac{6%
		}{q}}({\tau -\rho })^{-1} \\
	& +C\left( \frac{\sqrt{\ln R}}{\tau -\rho }\right) ^{\frac{2}{q}}+C\left(
	\frac{\sqrt{\ln R}}{\tau -\rho }\right) \tau ^{2-\frac{4}{q}}+C\left( \frac{R%
	}{\tau -\rho }\right) ^{%
		\frac{2q^{\prime }}{q}}\left(\Vert \nabla v\Vert _{L^{q}(B_{R}\backslash B_{R/2})}\right) ^{%
		\frac{2q^{\prime }}{q}+2}.
\end{align*}%
Then an application of Lemma \ref{Gia iteration} yields
\begin{equation*}
	\int_{B_{R/2}}|\nabla w|^{2}dx\leq CR^{-\frac{4}{q}}+C\left( \frac{\sqrt{%
			\ln R}}{R}\right) +C\left( \Vert \nabla v\Vert _{L^{q}(B_{R}\backslash
		B_{R/2})}\right) ^{\frac{2q^{\prime }}{q}+2}.
\end{equation*}%
Letting $R\rightarrow \infty $, by noting that (\ref{eq:nabla v}) we
have
\begin{equation*}
	\nabla w\equiv 0,
\end{equation*}%
and $\tilde{w}\equiv 0$. Similar arguments as in {\bf Step 1}, we complete the
proof.

%\smallskip
%
%\noindent{\textbf{Proof of Case II}}:
%
%\smallskip

\section{Proof of Theorem \ref{th3}}
\smallskip

%\textbf{Proof of Case IV}

%In this case, we want to prove that $\nabla v\in L^{2}(\Omega ),$ then the
%proof is complete by \textbf{Case II}.

In this case, we want to prove that $\nabla v\in L^2(\Omega),$ then the proof is complete by Theorem \ref{th2}.

%{\bf Proof of Theorem \ref{th3}.}

\bigskip
{\bf \underline{Case of $v\in L^p$ with  $1< p\leq 2$}}.
\bigskip

 It's similar to Theorem \ref{th1}. In fact, we let $\phi(x)\in C_0^\infty(B_R)$ with $0\leq \phi\leq 1$ as in (\ref{eq:phi}). Since
\beno
\int_{B_\tau\setminus B_1}\nabla (\phi v)dx=-\int_{\partial B_1}n\cdot v\phi dS=-\phi(1)\int_{\partial B_1}n\cdot v dS=0,
\eeno
due to (\ref{BC-slip}),
one could take a vector-valued function $\hat{w}: B_\tau\setminus{B_{\frac23\tau}}\rightarrow \mathbb{R}^2$ such that
$\hat{w}\in W^{1,s}_0(B_\tau\setminus{B_{\frac23\tau}})$ and
$\nabla\cdot \hat{w}=\nabla\cdot[\phi {v}]$ as in (\ref{esti-ws'}). Then
making the inner products $(\phi {v}-\hat{w})$ on both sides of the equation (\ref{NS-2D-r}), by $\nabla\cdot \hat{w}=\nabla\cdot[\phi{v}]$ we have
\beno
&&\hspace{-2mm}\int_{B_\tau\setminus{\overline{B_1}}}\phi|\nabla v|^2\,dx
\\&\leq& -\int_{B_\tau\setminus{\overline{B_1}}}\nabla\phi\cdot\nabla v \cdot {v} \,dx+\int_{B_\tau\setminus{\overline{B_1}}}\nabla \hat{w}:\nabla v  \,dx
-\int_{B_\tau\setminus{\overline{B_1}}}u\cdot\nabla v \cdot \phi {v} \,dx\\
&& +\int_{B_\tau\setminus{\overline{B_1}}}u\cdot\nabla v \cdot \hat{w} \,dx-\int_{B_\tau\setminus{\overline{B_1}}}v\cdot\nabla u_0 \cdot \phi {v} \,dx+\int_{B_\tau\setminus{\overline{B_1}}}v\cdot\nabla u_0 \cdot \hat{w} \,dx\\ && -\int_{\partial B_1}\frac{\partial v}{\partial n}\cdot v \phi dx
\doteq  I_1+\cdots+I_7,
\eeno
where $I_7\leq C$, the boundary value of the velocity at $\partial B_1$ is bounded. Similarly,
\beno
I_3&\leq& \frac12\left|\int_{B_\tau\setminus{\overline{B_1}}}u\cdot\nabla\phi |v|^2 \,dx\right|+C
\leq C(\tau-\rho)^{-1}\int_{B_{\tau }\setminus {B_{\frac{2}{3}\tau }}}|v|^3dx+C
\eeno
and other terms are similar to  Theorem \ref{th1}. The same arguments yield that
\ben\label{eq:L2 estimate'}
\int_{B_{R/2}\backslash{\overline{B_1}}}|\nabla v|^2\,dx\leq C.
\een
Applying Theorem \ref{th2}, we complete the proof.
\bigskip
\underline{{\bf  Case of $v\in L^p$ with $2<p<\infty$}}.
\bigskip

Let $\phi (x)\in C_{0}^{\infty }(\mathbb{R}^{2})$ be a
cut-off function defined as in (\ref{eq:phi}). Write $\tilde{w}^{2q}=(%
\tilde{w}^{2})^{q}$. For $q\geq 1$, by H\"{o}lder and Young inequalities we have
\ben\label{eq:wq}
	\int_{\Omega }\tilde{w}^{2q}\phi ^{2q }dx &=&\int_{\Omega
	}(v_{2},-v_{1})\cdot \nabla \lbrack \tilde{w}^{2q-2}\tilde{w}\phi ^{2q }]dx \nonumber\\
	&\leq &(2q-1)\int_{\Omega }|v||\nabla \tilde{w}|\tilde{w}^{2q-2}\phi ^{2q
	}dx+2q \int_{\Omega }|v||\nabla \phi ||\tilde{w}|^{2q-1}\phi ^{2q -1}dx
	\nonumber\\
	&\leq &\frac{1}{2}\int_{\Omega }\tilde{w}^{2q}\phi ^{2q
	}dx+C(q)\|v\|_{2q}^{\frac{2q}{q+1}} \left(\int_{\Omega }|\nabla \tilde{w}|^{2}\tilde{w}^{2q-2}\phi
	^{2q }dx\right)^{\frac{q}{q+1}} \nonumber\\
&&+ C(q)\|v\|_{2q}^{2q} (\tau-\rho)^{-2q}.
\een
%Due to the growth estimates $v\in \chi ^{\alpha ,\infty }$, we have
%\begin{align}
%	\int_{B_{1}^{c}}\tilde{w}^{2q}\phi ^{2q }dx& \leq C(q)R^{2\alpha
%	}\int_{\Omega }|\nabla \tilde{w}|^{2}\tilde{w}^{2q-4}\phi ^{2q }dx  \notag
%	\\
%	& +C(q )R^{\alpha -1}\int_{\Omega }|\tilde{w}|^{2q-1}\phi ^{2q -1}dx.
%	\label{w estimate-mhd4}
%\end{align}%
On the other hand, multiply $\phi ^{2q }\tilde{w}^{2q-2}\tilde{w}$ on both sides
of (\ref{NSv-2D}), and we have
\begin{align}
	& (2q-1)\int_{\Omega }|\nabla \tilde{w}|^{2}\tilde{w}^{2q-2}\phi ^{2q }dx
	\notag \\
	& \leq \frac{1}{2q}\int_{\Omega }\tilde{w}^{2q}\Delta (\phi ^{2q })dx+%
	\frac{1}{2q}\int_{\Omega }\tilde{w}^{2q}v\cdot \nabla (\phi ^{2q })dx
	\notag \\
	& +\frac{1}{2q}\int_{\Omega }\tilde{w}^{2q}u_{0}\cdot \nabla (\phi
	^{2q })dx\notag\\
	& \doteq II_{1}+\cdots +II_{3},  \label{energy estimate-w}
\end{align}%
and the last term vanishes. For the first two terms, there hold
\beno
II_1\leq C (\tau-\rho)^{-2}\int_{B_{\tau
		}\backslash {B_{1}}}|\tilde{w}|^{2q} dx,
\eeno
and
\beno
II_2\leq C (\tau-\rho)^{-1}\|v\|_{2q}\|\tilde{w}\|_{L^{\frac{4q^2}{2q-1}}(B_{\tau }\setminus {B_{\frac{2}{3}\tau }} )}^{2q}.
\eeno
Noting that
\beno
\|\tilde{w}^{q}\|_{L^{\frac{4q}{2q-1}}(B_{\tau }\setminus {B_{\frac{2}{3}\tau }} )}\leq C \|\nabla(\tilde{w}^{q})\|_{L^{2}(B_{\tau }\setminus {B_{\frac{2}{3}\tau }} )}^{\frac1{2q}}\|\tilde{w}^{q}\|_{L^{2}(B_{\tau }\setminus {B_{\frac{2}{3}\tau }} )}^{1-\frac1{2q}}+C\tau^{-\frac1{2q}}\|\tilde{w}^{q}\|_{L^{2}(B_{\tau }\setminus {B_{\frac{2}{3}\tau }} )},
\eeno
which implies
\beno
II_2&\leq& C (\tau-\rho)^{-1}\|v\|_{2q}\|\nabla(\tilde{w}^{q})\|_{L^{2}(B_{\tau }\setminus {B_{\frac{2}{3}\tau }} )}^{\frac1{q}}\|\tilde{w}^{q}\|_{L^{2}(B_{\tau }\setminus {B_{\frac{2}{3}\tau }} )}^{2-\frac1{q}}\\
&&+C (\tau-\rho)^{-1}\tau^{-\frac1{q}}\|v\|_{2q}\|\tilde{w}^{q}\|_{L^{2}(B_{\tau }\setminus {B_{\frac{2}{3}\tau }} )}^2.
\eeno
Then
\ben\label{energy estimate-nablaw}
	&&\int_{\Omega }|\nabla \tilde{w}|^{2}\tilde{w}^{2q-2}\phi ^{2q }dx
	\notag \\
	& \leq& C (\tau-\rho)^{-2}\int_{B_{\tau
		}\backslash {B_{1}}}|\tilde{w}|^{2q} dx \nonumber\\
&&+  C (\tau-\rho)^{-1}\|v\|_{2q}\|\nabla(\tilde{w}^{q})\|_{L^{2}(B_{\tau }\setminus {B_{\frac{2}{3}\tau }} )}^{\frac1{q}}\|\tilde{w}^{q}\|_{L^{2}(B_{\tau }\setminus {B_{\frac{2}{3}\tau }} )}^{2-\frac1{q}}\nonumber\\
&&+C (\tau-\rho)^{-1}\tau^{-\frac1{q}}\|v\|_{2q}\|\tilde{w}^{q}\|_{L^{2}(B_{\tau }\setminus {B_{\frac{2}{3}\tau }} )}^2.
\een
Substituting these estimates into the previous estimate of (\ref{eq:wq}), we get
\ben\label{eq:wq2}
	\int_{\Omega }\tilde{w}^{2q}\phi ^{2q }dx
	&\leq &C(q)\|v\|_{2q}^{\frac{2q}{q+1}} \left(  (\tau-\rho)^{-2}\int_{B_{\tau
		}\backslash {B_{1}}}|\tilde{w}|^{2q} dx \right)^{\frac{q}{q+1}} \nonumber\\
&&+ C(q)\|v\|_{2q}^{\frac{2q}{q+1}} \left(  (\tau-\rho)^{-1}\|v\|_{2q}\|\nabla(\tilde{w}^{q})\|_{L^{2}(B_{\tau }\setminus {B_{\frac{2}{3}\tau }} )}^{\frac1{q}}\|\tilde{w}^{q}\|_{L^{2}(B_{\tau }\setminus {B_{\frac{2}{3}\tau }} )}^{2-\frac1{q}}\right)^{\frac{q}{q+1}}\nonumber\\
&&+ C(q)\|v\|_{2q}^{\frac{2q}{q+1}} \left(  (\tau-\rho)^{-1}\tau^{-\frac1{q}}\|v\|_{2q}\|\tilde{w}^{q}\|_{L^{2}(B_{\tau }\setminus {B_{\frac{2}{3}\tau }} )}^2 \right)^{\frac{q}{q+1}} \nonumber\\
&&+ C(q)\|v\|_{2q}^{2q} (\tau-\rho)^{-2q}.
\een
Combining (\ref{energy estimate-nablaw}) and \ref{eq:wq2}, let $2q=p$ and the known condition $\|v\|_{2q}\leq C$ implies that
\beno
&&\int_{B_{\rho
		}\backslash {B_{1}}}|\tilde{w}|^{2q} dx+\left(\int_{B_{\rho
		}\backslash {B_{1}}}|\nabla(\tilde{w}^q)|^{2} dx\right)^{\frac{2q+1}{2q+2}}\nonumber\\
&\leq &\frac12 \int_{B_{\tau
		}\backslash {B_{1}}}|\tilde{w}|^{2q} dx+\frac12\left(\int_{B_{\tau
		}\backslash {B_{1}}}|\nabla(\tilde{w}^q)|^{2} dx\right)^{\frac{2q+1}{2q+2}}+ C(q)(\tau-\rho)^{-2q}+ C(q)(\tau-\rho)^{-q} \nonumber\\
&&+ C(q)(\tau-\rho)^{-4q}+  C(q)(\tau-\rho)^{-2q\frac{2q+1}{4q+1}}+ C(q)(\tau-\rho)^{-2q\frac{2q+1}{2q-1}},\nonumber\\
\eeno
where we used $\tau>1$ and Young's inequality with the index
\beno
\frac{4q+1}{2(q+1)(2q+1)}+  \frac{2q-1}{2(q+1)}+\frac{1}{2q+1}=1,
\eeno
and
\beno
\frac{2q-1}{4q(q+1)}+\frac{1}{2q}+  \frac{(2q-1)(2q+1)}{4q(q+1)}=1.
\eeno
Applying Lemma \ref{Gia iteration} again, we have
\begin{equation*}
	\tilde{w}=0.
\end{equation*}%
Recall $v'$ and $w'$ and estimate it as in (\ref{eq:w'}), then it follows that $\nabla v\in L^p$ for any $p>1$ due to $\Delta v'=\nabla ^{\perp }w'$.
Applying  Theorem \ref{th2}, the proof is complete.

\section{Proof of Theorem \protect\ref{th4}}

\label{proof1}

\smallskip
{\bf {\underline {Case 1: $u(x)$ is a $q$-generalized solution for $1<q\leq 3/2.$}}}
\bigskip

It is obvious to see that the pressure $\pi$ is almost silent in the definition
of $q$-generalized solutions, which conceals some information about
pressure.  Note that in exterior domains, the Calder\'{o}%
n-Zygmund inequality does not work, and the estimate for pressure via velocity is not available. One can apply the techniques used in Theorem \ref{th1} to prove the triviality of the $q$-generalized solutions. Similar as the local energy estimate (\ref{eq:energy}), one can get
\begin{align}\label{eq:energy}
&\int_{B_\tau\setminus{\overline{B_1}}}\phi|\nabla u|^2\,dx \notag \\
&\leq\frac14\|\nabla u\|_{L^2(B_\tau\setminus{\overline{B_1}})}^2+ \frac{C}{(\tau-\rho)^{2}}\|u\|^{2}_{L^{2}(B_\tau\setminus{B_{\frac23\tau}})}+\frac{C}{\tau-\rho}\|u\|^{3}_{L^{3}(B_\tau\setminus{B_{\frac23\tau}})},
\end{align}
which implies the required result by using the same arguments, since $u\in L^p$ with $2<p\leq 6$.
%As a matter of fact, the current case is a particular one that was stated in Theorem \ref{th2}, where we just assume that $a=0$ and $C=0$ in the boundary conditions.
So, for concision, we skip the details.

\bigskip
{\bf \underline{Case 2: $u\mathbf{(}x\mathbf{)}\in BMO^{-1}(\Omega )$ for $3/2<q<2.$}}
\bigskip

We construct a
cut-off radially nonincreasing function $\zeta _{R}(x)\in C_{0}^{\infty }(%
\mathbb{R}^{2})$ for $R\gg 1$ by $0\leq \zeta _{R}(x)\leq 1$ which satisfies the
followings%
\begin{equation*}
	\zeta _{R}(x)=\left\{
	\begin{array}{c}
		1,\text{ \ }x\in B_{\rho } \\
		0,\text{ \ }x\in B_{\tau }^{c}%
	\end{array}%
	\right. ,
\end{equation*}%
with
\begin{equation*}
	\frac{R}{2}<\frac{\tau }{2}\leq R<\rho <\tau <2R,
\end{equation*}%
moreover,
\begin{equation*}
	\left\Vert \nabla \zeta _{R}(x)\right\Vert _{L^{\infty }(\mathbb{R}%
		^{2})}\leq \frac{C}{\tau -\rho },\text{ \ }\left\Vert \nabla ^{2}\zeta
	_{R}(x)\right\Vert _{L^{\infty }(\mathbb{R}^{2})}\leq \frac{C}{\left( \tau
		-\rho \right) ^{2}},
\end{equation*}%
where $C$ is independent of $x$ and $R.$ It is easy to know that $\nabla
\zeta _{R}(x)$ is supported in $\mathcal{A}_{R}=B_{\tau }\backslash \overline{B}%
_{\rho }.$
Now, multiplying both sides of (\ref{EQ:NS}) with $\zeta _{R}u-\varphi $, and
noticing that $\zeta _{R}u-\varphi $ is divergence-free$,$ one has%
\begin{eqnarray}
\int_{B_{\tau }\backslash \bar{B}_{1}}|\nabla u|^{2}\zeta _{R}dx
&=&-\int_{B_{\tau }\backslash B_{\tau /4}}\nabla u\cdot \nabla \zeta
_{R}\cdot udx+\int_{B_{\tau }\backslash B_{\tau /4}}\nabla u:\nabla \varphi
dx  \notag \\
&&-\int_{B_{\tau }\backslash B_{\tau /4}}\left( u\cdot \nabla \right) u\cdot
\zeta _{R}udx+\int_{B_{\tau }\backslash B_{\tau /4}}\left( u\cdot \nabla
\right) u\cdot \varphi dx  \notag \\
&=&\frac{1}{2}\int_{B_{\tau }\backslash B_{\tau /4}}|u|^{2}\Delta \zeta
_{R}dx+\int_{B_{\tau }\backslash B_{\tau /4}}\nabla u:\nabla \varphi dx
\notag \\
&&+\frac{1}{2}\int_{B_{\tau }\backslash B_{\tau /4}}|u|^{2}(u\cdot \nabla
\zeta _{R})dx+\int_{B_{\tau }\backslash B_{\tau /4}}\left( u\cdot \nabla
\right) u\cdot \varphi dx  \notag \\
&:=&J_{1}+J_{2}+J_{3}+J_{4}.  \label{J1-4}
\end{eqnarray}%
We aim to prove that each $J_{j}$ $(j=1,2,3,4)$ tends to zero as $R$ goes to
infinity, which implies that $\left\Vert \nabla u\right\Vert _{L^{2}(\Omega
)}=0,$ then it implies that $u\equiv 0.$ Firstly, it follows by H\"{o}lder's
inequality that%
\begin{eqnarray*}
J_{1} &\leq &C\left\Vert u\right\Vert _{L^{\frac{2q}{2-q}}\left( B_{\tau
}\backslash B_{\tau /4}\right) }^{2}\left\Vert \Delta \zeta _{R}\right\Vert
_{L^{\frac{q}{2q-2}}\left( B_{\tau }\backslash B_{\tau /4}\right) }\\ &\leq&
\frac{C}{(\tau -\rho )^{2}}\left\Vert u\right\Vert _{L^{\frac{2q}{2-q}%
}\left( B_{\tau }\backslash B_{\tau /4}\right) }^{2}(\tau -\rho )^{\left( 4-%
\frac{4}{q}\right) }
\leq \frac{C}{R^{\left( \frac{4}{q}-2\right) }}\left\Vert u\right\Vert
_{L^{\frac{2q}{2-q}}\left( B_{\tau }\backslash B_{\tau /4}\right) }^{2},
\end{eqnarray*}%
and%
\begin{align*}
J_{2}& \leq \left\Vert \nabla u\right\Vert _{L^{q}(B_{\tau }\backslash \bar{B%
}_{1})}\left\Vert \nabla \varphi \right\Vert _{L^{\frac{q}{q-1}}(B_{\tau
}\backslash B_{\tau /4})}\leq \left\Vert \nabla u\right\Vert _{L^{q}(B_{\tau
}\backslash \bar{B}_{1})}\left\Vert \nabla \zeta _{R}\cdot u\right\Vert _{L^{%
\frac{q}{q-1}}(B_{\tau }\backslash B_{\tau /4})} \\
& \leq \frac{C}{(\tau -\rho )}\left\Vert \nabla u\right\Vert _{L^{q}(B_{\tau
}\backslash \bar{B}_{1})}\left\Vert u\right\Vert _{L^{\frac{2q}{2-q}%
}(B_{\tau }\backslash B_{\tau /4})}(\tau -\rho )^{\left( 3-\frac{4}{q}%
\right) } \\
& \leq \frac{C}{R^{\left( \frac{4}{q}-2\right) }}\left\Vert \nabla
u\right\Vert _{L^{q}(B_{\tau }\backslash \bar{B}_{1})}\left\Vert
u\right\Vert _{L^{\frac{2q}{2-q}}(B_{\tau }\backslash B_{\tau /4})}.
\end{align*}%
The estimate of $J_{3}$ is slightly different, we proceed with the
assumption that $u\in BMO^{-1}(\Omega ).$ Since $u\in BMO^{-1}(\Omega ),$
then each component of $u$ can be represented by%
\begin{equation*}
u_{i}=\sum\limits_{j=1}^{2}\partial _{j}g_{j}^{i},\text{ \ }i=1,2,
\end{equation*}%
for some suitable functions $g_{j}^{i}\in BMO(\Omega ).$ The estimate of $%
J_{3}$ is given as follows.
\begin{align*}
J_{3}& =\frac{1}{2}\int_{B_{\tau }\backslash B_{\tau /4}}|u|^{2}(u\cdot
\nabla \zeta _{R})dx=\frac{1}{2}\sum\limits_{i=1}^{2}\int_{B_{\tau
}\backslash B_{\tau /4}}|u|^{2}(u_{i}\partial _{i}\zeta _{R})dx \\
& =\frac{1}{2}\sum\limits_{1\leq i,j\leq 2}\int_{B_{\tau }\backslash B_{\tau
/4}}|u|^{2}\partial _{j}(g_{j}^{i}-[g_{j}^{i}]_{\tau })\partial _{i}\zeta
_{R}dx\\
&=-\frac{1}{2}\sum\limits_{1\leq i,j\leq 2}\int_{B_{\tau }\backslash
B_{\tau /4}}\partial _{j}\left( |u|^{2}\partial _{i}\zeta _{R}\right)
(g_{j}^{i}-[g_{j}^{i}]_{\tau })dx \\
& =-\frac{1}{2}\sum\limits_{1\leq i,j\leq 2}\int_{B_{\tau }\backslash
B_{\tau /4}}|u|^{2}\partial _{ij}^{2}\zeta _{R}(g_{j}^{i}-[g_{j}^{i}]_{\tau
})dx\\
&\hspace{5mm}-\sum\limits_{1\leq i,j\leq 2}\int_{B_{\tau }\backslash B_{\tau
/4}}\left( u\cdot \partial _{j}u\right) \partial _{i}\zeta
_{R}(g_{j}^{i}-[g_{j}^{i}]_{\tau })dx \\
& :=J_{31}+J_{32},
\end{align*}%
where $[g_{j}^{i}]_{\tau }$ is the mean value of $\displaystyle\int_{B_{\tau
}\backslash B_{\tau /4}}g_{j}^{i}dx$ on $B_{\tau }\backslash B_{\tau /4}.$
Then%
\begin{align*}
|J_{31}|& \leq \sup_{i,j}\left\Vert g_{j}^{i}-[g_{j}^{i}]_{\tau }\right\Vert
_{L^{\frac{2q}{3q-3}}(B_{\tau }\backslash B_{\tau /4})}\left\Vert
u\right\Vert _{L^{\frac{2q}{2-q}}(B_{\tau }\backslash B_{\tau
/4})}^{2}\left\Vert \Delta \zeta _{R}\right\Vert _{L^{\frac{2q}{q-1}%
}(B_{\tau }\backslash B_{\tau /4})} \\
& \leq CR^{2-\frac{4}{q}}\sup_{i,j}\left\Vert g_{j}^{i}\right\Vert
_{BMO(B_{\tau }\backslash B_{\tau /4})}\left\Vert u\right\Vert _{L^{\frac{2q%
}{2-q}}(B_{\tau }\backslash B_{\tau /4})}^{2},
\end{align*}%
where we used the inequality (\ref{BMO pro}), and similarly%
\begin{align*}
|J_{32}|& \leq \sup_{i,j}\left\Vert g_{j}^{i}-[g_{j}^{i}]_{\tau }\right\Vert
_{L^{\frac{2q}{2q-3}}(B_{\tau }\backslash B_{\tau /4})}\left\Vert
u\right\Vert _{L^{\frac{2q}{2-q}}(B_{\tau }\backslash B_{\tau
/4})}\left\Vert \nabla u\right\Vert _{L^{q}(B_{\tau }\backslash B\rho
)}\left\Vert \nabla \zeta _{R}\right\Vert _{L^{\frac{2q}{q-1}}(B_{\tau
}\backslash B_{\tau /4})} \\
& \leq CR^{2-\frac{4}{q}}\sup_{i,j}\left\Vert g_{j}^{i}\right\Vert
_{BMO(B_{\tau }\backslash B_{\tau /4})}\left\Vert u\right\Vert _{L^{\frac{2q%
}{2-q}}(B_{\tau }\backslash B_{\tau /4})}\left\Vert \nabla u\right\Vert
_{L^{q}(B_{\tau }\backslash B_{\tau /4})}.
\end{align*}%
Note that $3/2<q<2,$ then%
\begin{equation*}
2q-3>0,\text{ }2-q>0,\text{ }2-\frac{4}{q}\leq 0.
\end{equation*}%
Finally, we estimate $J_{4}.$
\begin{align*}
J_{4}& =\int_{B_{\tau }\backslash B_{\tau /4}}\left( u\cdot \nabla \right)
u\cdot \varphi dx=\sum\limits_{1\leq i,j\leq 2}\int_{B_{\tau }\backslash
B_{\tau /4}}u_{i}\partial _{i}u_{j}\varphi _{j}dx \\
& =-\sum\limits_{1\leq i,j,l\leq 2}\int_{B_{\tau }\backslash B_{\tau
/4}}\partial _{l}(g_{l}^{i}-[g_{l}^{i}]_{\tau })u_{j}\partial _{i}\varphi
_{j}dx \\
& =\sum\limits_{1\leq i,j,l\leq 2}\int_{B_{\tau }\backslash B_{\tau
/4}}(g_{l}^{i}-[g_{l}^{i}]_{\tau })\partial _{l}\left( u_{j}\partial
_{i}\varphi _{j}\right) dx \\
& =\sum\limits_{1\leq i,j,l\leq 2}\int_{B_{\tau }\backslash B_{\tau
/4}}(g_{l}^{i}-[g_{l}^{i}]_{\tau })\left( \partial _{l}u_{j}\partial
_{i}\varphi _{j}+u_{j}\partial _{il}^{2}\varphi _{j}\right) dx \\
& :=J_{41}+J_{42}.
\end{align*}%
In a similar way of $J_{3,}$ one has%
\begin{align*}
J_{41}& \leq \frac{C}{R}\sup_{i,j}\left\Vert
g_{j}^{i}-[g_{j}^{i}]_{R}\right\Vert _{L^{\frac{2q}{3q-4}}(B_{2R})}\left%
\Vert \nabla u\right\Vert _{L^{q}(B_{\tau }\backslash \bar{B}%
_{1})}\left\Vert u\right\Vert _{L^{\frac{2q}{2-q}}(B_{\tau }\backslash B\rho
)} \\
& \leq CR^{2-\frac{4}{q}}\sup_{i,j}\left\Vert g_{j}^{i}\right\Vert
_{BMO(B_{2R})}\left\Vert u\right\Vert _{L^{\frac{2q}{2-q}}(B_{\tau
}\backslash B\rho )}\left\Vert \nabla u\right\Vert _{L^{q}(B_{\tau
}\backslash B\rho )},
\end{align*}%
and%
\begin{align*}
J_{42} \leq &\frac{C}{R^{2}}\sup_{i,j}\left\Vert g_{j}^{i}-[g_{j}^{i}]_{\tau
}\right\Vert _{L^{\frac{q}{2q-2}}(B_{\tau }\backslash B_{\tau
/4})}\left\Vert u\right\Vert _{L^{\frac{2q}{2-q}}(B_{\tau }\backslash
B_{\tau /4})}^{2} \\
&+C\sup_{i,j}\left\Vert g_{j}^{i}-[g_{j}^{i}]_{\tau }\right\Vert _{L^{\frac{%
2q}{3q-4}}}\left\Vert u\right\Vert _{L^{\frac{2q}{2-q}}(B_{\tau }\backslash
B_{\tau /4})}\left\Vert \nabla \zeta _{R}\cdot \nabla u\right\Vert
_{_{L^{q}(B_{\tau }\backslash B_{\tau /4})}} \\
\leq &CR^{2-\frac{4}{q}}\sup_{i,j}\left\Vert g_{j}^{i}\right\Vert
_{BMO(B_{\tau }\backslash B_{\tau /4})}\left\Vert u\right\Vert _{L^{\frac{2q%
}{2-q}}(B_{\tau }\backslash B_{\tau /4})}^{2} \\
&+CR^{2-\frac{4}{q}}\sup_{i,j}\left\Vert g_{j}^{i}\right\Vert _{BMO(B_{\tau
}\backslash B_{\tau /4})}\left\Vert u\right\Vert _{L^{\frac{2q}{2-q}%
}(B_{\tau }\backslash B_{\tau /4})}\left\Vert \nabla u\right\Vert
_{L^{q}(B_{\tau }\backslash B_{\tau /4})}.
\end{align*}%
It is obvious to see that $J_{1},$ $J_{2},$ $J_{3}$ and $J_{4}$ tend to zero
as $R$ goes to $\infty ,$ so we conclude that $u\equiv 0.$

\bigskip
{\bf{\underline{Case 3: $u(x)\in L^{4}(\Omega)$ for a $D$-solution}}}
\bigskip

As discussed in Case 1, the estimate of pressure is important to obtain suitable estimates in this
situation. To this end, multiplying the Navier-Stokes equations (\ref{EQ:NS})
by $\tilde{\psi}\in C_{0}^{\infty }(\Omega )$ (not necessarily solenoidal),
then integrating by parts yields
\begin{equation}
	(\nabla u,\nabla \tilde{\psi})=-(u\cdot \nabla u,\tilde{\psi})+(\pi,\nabla
	\cdot \tilde{\psi}).  \label{pressure}
\end{equation}%
If the convective term $u\cdot \nabla u$ has a mild degree of regularity, to
every $q$-generalized solution we are able to associate a pressure $\pi$ such
that (\ref{pressure}) holds. Therefore, for a locally Lipschitz exterior
domain of $\mathbb{R}^{2}$, if $u\cdot \nabla u\in D_{0}^{-1,q}(\Omega ),$
there exists a unique $\pi\in L^{q}(\Omega )$ satisfying (\ref{pressure}) for
all $\tilde{\psi}\in C_{0}^{\infty }(\Omega ).$ Furthermore, the following
inequality holds
\begin{equation}
	\left\Vert \pi\right\Vert _{L^{q}}\leq C\left( \left\vert u\cdot \nabla
	u\right\vert _{-1,q}+\left\Vert \nabla u\right\Vert _{L^{q}}\right) .
	\label{d-1}
\end{equation}%
The proof of this argument can be found in \cite[Lemma V.1.1 in Page 305]%
{galdibook}. Now Let $\psi (\xi )$ be a nonincreasing smooth function
defined in $\mathbb{R}^{2}$ with $\psi (\xi )=1$ if $|\xi |\leq 1/2$ and $%
\psi (\xi )=0$ if $|\xi |\geq 1,$ and set, for $R$ large enough,
\begin{equation*}
	\psi _{R}(x)=\psi \left( \frac{\ln \ln |x|}{\ln \ln R}\right) ,\text{ \ }%
	x\in \Omega .
\end{equation*}%
Note that, for a suitable constant $c$ independent of $R,$ there holds
\begin{equation*}
	|\nabla \psi _{R}(x)|\leq \frac{c}{\ln \ln R}\frac{1}{|x|\ln |x|},
\end{equation*}%
and $\nabla \psi _{R}(x)\not\equiv 0,$ only if $x\in \tilde{\Omega}_{R},$
where
\begin{equation*}
	\tilde{\Omega}_{R}=\{x\in \Omega :\exp \sqrt{\ln R}<|x|<R\}.
\end{equation*}%
Multiplying equation (\ref{EQ:NS}) by $\psi _{R}u\mathbf{,}$ integrating by
parts over $\Omega $ and taking the divergence free condition into account
yield
\begin{equation*}
	\int_{\Omega }\psi _{R}u\cdot (u\cdot \nabla )u\mathbf{=-}\frac{1}{2}%
	\int_{\Omega }|u\mathbf{|}^{2}u\cdot \nabla \psi _{R}dx,
\end{equation*}%
\begin{equation*}
	\int_{\Omega }\psi _{R}u\cdot \nabla \pi=\mathbf{-}\int_{\Omega }\pi\nabla \psi
	_{R}\cdot udx,
\end{equation*}%
and
\begin{equation*}
	\int_{\Omega }\psi _{R}u\cdot \Delta udx\mathbf{=-}\int_{\Omega }\nabla \psi
	_{R}\cdot \nabla u\cdot udx-\int_{\Omega }\psi _{R}\nabla u:\nabla udx%
	\mathbf{,}
\end{equation*}%
where we have used the zero boundary condition (\ref{boundary}). Then it
follows that
\begin{equation}
	\int_{\Omega }\psi _{R}\nabla u:\nabla u=\mathbf{-}\int_{\Omega }\nabla \psi
	_{R}\cdot \nabla u\cdot udx+\int_{\Omega }\pi\nabla \psi _{R}\cdot udx+\frac{1%
	}{2}\int_{\Omega }|u\mathbf{|}^{2}u\cdot \nabla \psi _{R}dx.
	\label{gradient}
\end{equation}%
By H\"{o}lder's inequality, one has
\begin{align}
	\int_{B_{R}\backslash \bar{B}_{1}}\psi _{R}\left\vert \nabla u\right\vert
	^{2}dx& =\int_{\Omega }\psi _{R}\left\vert \nabla u\right\vert ^{2}dx\leq
	||\nabla \psi _{R}u||_{L^{2}(\tilde{\Omega}_{R})}||\nabla u||_{L^{2}(\Omega
		)}  \notag \\
	& +||\pi||_{L^{2}(\Omega )}||\nabla \psi _{R}u||_{L^{2}(\tilde{\Omega}%
		_{R})}+C||u||_{L^{4}(\Omega )}^{2}||\nabla \psi _{R}u||_{L^{2}(\tilde{\Omega}%
		_{R})},  \label{ineqs}
\end{align}%
We found by definition that
\begin{equation}
	\left\Vert u\cdot \nabla u\right\Vert _{-1,2}=\sup_{w\in D_{0}^{1,2}(\Omega
		);||\nabla w||_{L^{2}}=1}\int_{\Omega }u\cdot \nabla u\cdot w\leq \left\Vert
	u\right\Vert _{L^{4}}^{2}\left\Vert \nabla w\right\Vert _{L^{2}}.
	\label{defi}
\end{equation}%
Since we assume that $u\in L^{4},$ (\ref{defi}) implies that $u\cdot \nabla
u\in D_{0}^{-1,2}(\Omega ).$ Thus, by (\ref{d-1}) one has
\begin{equation}
	\pi\in L^{2}(\Omega ).  \label{pest}
\end{equation}%
On the other hand, it follows that

\begin{align}
	||\nabla \psi _{R}u||_{L^{2}(\tilde{\Omega}_{R})}^{2}& \leq \frac{c_{1}}{%
		\left( \ln \ln R\right) ^{2}}\int_{\exp \sqrt{\ln R}}^{R}\frac{|u|^{2}}{%
		\left( |z|\ln |z|\right) ^{2}}dz  \notag \\
	& =\frac{c_{1}}{\left( \ln \ln R\right) ^{2}}\int_{\exp \sqrt{\ln R}%
	}^{R}\int_{0}^{2\pi }\frac{|u(r,\theta )|^{2}r^{-1}}{\left( \ln |r|\right)
		^{2}}drd\theta .  \label{lq2}
\end{align}%
It remains to estimate the right-hand side of (\ref{lq2}). Note that
\begin{equation*}
	u(r,\theta )=u(r_{0},\theta )+\int_{r_{0}}^{r}\frac{\partial u}{\partial \xi
	}d\xi ,\text{ \ for }r\geq r_{0}>1.
\end{equation*}%
Then
\begin{align*}
	|u(r,\theta )|^{2}& =\left( u(r_{0},\theta )+\int_{r_{0}}^{r}\frac{\partial u%
	}{\partial \xi }d\xi \right) ^{2} \\
	& \leq 2\left( \left\vert u(r_{0},\theta )\right\vert ^{2}+\left\vert
	\int_{r_{0}}^{r}\frac{\partial u}{\partial \xi }d\xi \right\vert ^{2}\right)
	.
\end{align*}%
It follows that
\begin{align*}
	\int_{0}^{2\pi }|u(r,\theta )|^{2}d\theta & \leq 2\left( \int_{0}^{2\pi
	}\left\vert u(r_{0},\theta )\right\vert ^{2}d\theta +\int_{0}^{2\pi
	}\left\vert \int_{r_{0}}^{r}\frac{\partial u}{\partial \xi }d\xi \right\vert
	^{2}d\theta \right) \\
	& \leq 2\left( \int_{0}^{2\pi }\left\vert u(r_{0},\theta )\right\vert
	^{2}d\theta +\int_{0}^{2\pi }\left\vert \int_{r_{0}}^{r}\frac{\partial u}{%
		\partial \xi }d\xi \right\vert ^{2}d\theta \right) .
\end{align*}%
By H\"{o}lder's inequality
\begin{align*}
	\left\vert \int_{r_{0}}^{r}\frac{\partial u}{\partial \xi }d\xi \right\vert
	^{2}& =\left\vert \int_{r_{0}}^{r}\frac{\partial u}{\partial \xi }\xi ^{%
		\frac{1}{2}}\xi ^{-\frac{1}{2}}d\xi \right\vert ^{2} \leq \left(
	\int_{r_{0}}^{r}\left( \frac{\partial u}{\partial \xi }\right)
	^{2}\left\vert \xi \right\vert d\xi \right) \text{ln}r ,
\end{align*}%
Therefore, we obtain
\begin{equation*}
	\int_{0}^{2\pi }\left\vert \int_{r_{0}}^{r}\frac{\partial u}{\partial \xi }%
	d\xi \right\vert ^{2}d\theta \leq \int_{0}^{2\pi }G(r,\theta )d\theta \leq
	\text{ln}r\left\Vert \nabla u\right\Vert _{L^{2}(B_{r}\backslash
		B_{r_{0}})}^{2}.
\end{equation*}%
Then from (\ref{lq2}), one has
\begin{align*}
	||\nabla \psi _{R}u||_{L^{2}(\tilde{\Omega}_{R})}^{2}& \leq \frac{c_{2}}{%
		\left( \ln \ln R\right) ^{2}}\int_{\exp \sqrt{\ln R}}^{R}\frac{r^{-1}(\ln r+C)}{%
		\left( \ln r\right) ^{2}}dr  \notag \\
	& \leq \frac{c_{2}}{\ln \ln R},
\end{align*}%
which implies that
\begin{equation*}
	\lim_{R\rightarrow \infty }||\nabla \psi _{R}u||_{L^{2}(\tilde{\Omega}%
		_{R})}=0,
\end{equation*}%
and from inequality (\ref{ineqs}), we conclude that
\begin{equation}
	\lim_{R\rightarrow \infty }\int_{B_{R}\backslash \bar{B}_{1}}\psi
	_{R}\left\vert \nabla u\right\vert ^{2}dx=0.  \label{zeros}
\end{equation}%
Relations (\ref{zeros}) and (\ref{gradient}) imply, by the monotone
convergence theorem, $\nabla u\equiv 0.$ It follows that $u$ must
identically be $0$ by our assumption. Therefore, the
proof of Theorem \ref{th4} is complete.

\noindent {\bf Acknowledgments.}
Z. Guo was partially supported by Natural Science Foundation of Jiangsu Province under grant BK20201478 and Qinlan Project of Jiangsu Universities. W. Wang was supported by NSFC under grant 12071054,  National Support Program for Young Top-Notch Talents and by Dalian High-level Talent Innovation Project (Grant 2020RD09).

\end{document}